\documentclass{article}
\usepackage{amssymb,latexsym,amsmath}

\begin{document}

{\Large{
\begin{center}
{\bf A Simple Example of a New Class of Landen Transformations}
\end{center}

\medskip

\begin{center}
{\bf Dante Manna and Victor H. Moll}
\end{center}

}}

\newcommand{\nn}{\nonumber}
\newcommand{\ba}{\begin{eqnarray}}
\newcommand{\ea}{\end{eqnarray}}
\newcommand{\dz}{\frac{d}{dz}}
\newcommand{\E}{{\mathfrak{E}}}
\newcommand{\F}{{\mathfrak{F}}}
\newcommand{\Ro}{{\mathfrak{R}}}
\newcommand{\ift}{\int_{0}^{\infty}}
\newcommand{\ifft}{\int_{- \infty}^{\infty}}
\newcommand{\no}{\noindent}
\newcommand{\X}{{\mathbb{X}}}
\newcommand{\Q}{{\mathbb{Q}}}
\newcommand{\R}{{\mathbb{R}}}
\newcommand{\Y}{{\mathbb{Y}}}
\newcommand{\Ftwo}{{{_{2}F_{1}}}}
\newcommand{\realpart}{\mathop{\rm Re}\nolimits}
\newcommand{\imagpart}{\mathop{\rm Im}\nolimits}

\newtheorem{Definition}{\bf Definition}[section]
\newtheorem{Thm}[Definition]{\bf Theorem} 
\newtheorem{Example}[Definition]{\bf Example} 
\newtheorem{Lem}[Definition]{\bf Lemma} 
\newtheorem{Note}[Definition]{\bf Note} 
\newtheorem{Cor}[Definition]{\bf Corollary} 
\newtheorem{Prop}[Definition]{\bf Proposition} 
\newtheorem{Problem}[Definition]{\bf Problem} 

\medskip

\no
{\bf 1. INTRODUCTION}. The 
method of completing squares yields an elementary procedure to evaluate 
\begin{equation}
I  =  \ifft \frac{dx}{ax^{2}+bx + c}.
\label{1.1}
\end{equation}
\no
Write
\begin{equation}
ax^{2}+bx + c  =  a \left[ \left( x + \frac{b}{2a} \right)^{2} + 
\frac{4ac-b^{2}}{4a^{2}} \right],
\nn
\end{equation}
\no
and use a linear change of variables to obtain
\begin{equation}
\ifft \frac{dx}{ax^{2}+bx + c}  =  \frac{2}{\sqrt{4ac-b^{2}}} 
\ifft \frac{dx}{x^{2}+1} 
= \frac{2 \pi}{\sqrt{4ac-b^{2}}}. 
\label{integral2}
\end{equation}
\no
Observe that $4ac-b^{2}>0$ is required for the convergence of (\ref{1.1}).

The goal of this paper is to present a new proof of (\ref{integral2}). We
illustrate a technique that will apply to {\em any
rational integrand}. 
Providing new proofs of an elementary result, such as (\ref{integral2}), is 
usually an effective tool to 
introduce students to more interesting Mathematics. 
The method discussed here has a rich history that we describe in section 2.

It is an unfortunate fact that, despite our best efforts, evaluating 
definite integrals is not very much in fashion today. Thus we rephrase the 
previous evaluation as a question in dynamical systems: replace the parameters
$a, \, b,$ and $c$ in (\ref{1.1}) with new ones given by the rules
\ba
a_{n+1} & =  a_{n} \left[ \frac{(a_{n}+3c_{n})^{2} - 3b_{n}^{2} }
{(3a_{n}+c_{n})(a_{n}+3c_{n}) - b_{n}^{2}} \right],  
\label{system1}
 \\
b_{n+1} & =  b_{n} \left[ \frac{3(a_{n}-c_{n})^{2} - b_{n}^{2} }
{(3a_{n}+c_{n})(a_{n}+3c_{n}) - b_{n}^{2}} \right], \nn  \\
c_{n+1} & =  c_{n} \left[ \frac{(3a_{n}+c_{n})^{2} - 3b_{n}^{2} }
{(3a_{n}+c_{n})(a_{n}+3c_{n}) - b_{n}^{2}} \right], \nn
\ea
\no
with $a_{0}=a, \, b_{0}=b,$ and $c_{0}=c$. 
The reader is asked to check 
that (\ref{1.1}) is 
invariant under (\ref{system1}), that is,
\ba
\ifft \frac{dx}{a_{n+1}x^{2} + b_{n+1}x + c_{n+1}} & = & 
\ifft \frac{dx}{a_{n}x^{2} + b_{n}x + c_{n}}, 
\label{invariance1}
\ea
\no
and to prove that
\ba
\lim\limits_{n \to \infty} a_{n} = 
\lim\limits_{n \to \infty} c_{n} = \frac{1}{2} \sqrt{4ac-b^{2}},  \, \, \,
\lim\limits_{n \to \infty} b_{n} = 0. 
\label{limits1}
\ea

\no
Once this is done, we can pass to the limit in (\ref{invariance1}) and 
use the invariance of $I$ to obtain
\ba
I & = & \frac{\pi}{\lim\limits_{n \to \infty} a_{n}} 
= \frac{2 \pi}{\sqrt{4ac - b^{2} } }.
\label{limits2}
\ea
\no
This leads directly to a proof of (\ref{integral2}). The advantage of this 
method is that it generalizes to integrands of higher degree.

We call  (\ref{system1}) a 
{\em rational Landen transformation}. In section 2
we discuss the historical precedent and motivation
behind such transformations. This history 
connects (\ref{limits2}) to the magic of
the arithmetic-geometric mean, $\pi$, and a wonderful
numerical calculation of Gauss. 

The rest of the paper is devoted to a detailed proof of the Landen 
transformation: the invariance of the rational integral (\ref{invariance1}) 
and the evaluation of the limits in (\ref{limits1}). A 
scaling of the integrand that is a crucial step in producing
this transformation is presented in 
section 3. The following
section presents the trigonometrical aspects of this problem and completes 
the proof of (\ref{invariance1}). An algebraic 
calculation shows that the discriminant of the quadratic in (\ref{1.1})
is preserved; that is,
\ba
4ac - b^{2} & = & 4a_{1}c_{1} - b_{1}^{2}. 
\label{discriminant1}
\ea
\no
This invariance is used in 
section 5 to analyze the dynamics of
(\ref{system1}) and to establish (\ref{limits1}). 

\medskip

\no
{\bf 2. LANDEN TRANSFORMATIONS}.  Many 
of the evaluations encountered in integral calculus illustrate the fact 
that definite integrals correspond to special values of functions. For 
example, the last integral in (\ref{integral2}) is given by $\pi = 
\tan^{-1}(\infty) - \tan^{-1}(-\infty)$. Other special values appear in
elementary courses:
\ba
\int_{0}^{1} \frac{dx}{\sqrt{3-x^{2}}} & = & \sin^{-1} \left( \frac{1}
{\sqrt{3}} \right). 
\label{intsine}
\ea

The same is true for more complicated integrals. For instance, when
$0 < b < a <1$, 
\ba
G(a,b) := \int_{0}^{\pi/2} \frac{d \theta}
{ \sqrt{ a^{2} \cos^{2} \theta + b^{2} \sin^{2} \theta} } 
= \frac{1}{a} K(k), 
\label{elliptic1}
\ea
\no
with $k^{2} = 1 - b^{2}/a^{2}$. Here $K$ is the {\em complete elliptic integral
of the first kind} defined by
\ba
K(k) & = & \int_{0}^{1} \frac{dx}{\sqrt{(1-x^{2})(1-k^{2}x^{2})}}.
\label{elliptic2}
\ea
\no
Elliptic integrals appear at the center of classical analysis. Their
name comes from the fact 
that they provide explicit formulas for the length of an ellipse. 

The inverse of 
\ba
f(z)  & = & \int_{0}^{z} \frac{dx}{\sqrt{(1-x^{2})(1-k^{2}x^{2})}}
\nn
\ea
\no
is similar to $\sin z$, so (\ref{intsine}) and (\ref{elliptic2}) are not so
different after all. This new function is the {\em elliptic sine} (or
{\em sinus amplitudinus}) of Jacobi [15], denoted by $\text{ sn }z$. 
It completes the trilogy: $ \, \sin z$ (circular), $\text{ sinh }z$ 
(hyperbolic), and $\text{ sn }z$ (elliptic).  The question of evaluating 
definite integrals sometimes comes down to how many functions one knows.

Our complaint
that students today are exposed only to the most basic of functions is not 
new. Klein states [16, \, p. 294]: 
{\em When I was a student, abelian functions 
were, as an effect of the Jacobian tradition, considered the uncontested 
summit of mathematics and each of us was ambitious to make progress in this
field. And now? The younger generation hardly knows abelian 
functions.}\footnote{The authors learned of this quote 
from the preface of [5].}\\

Suppose that $a$ and $b$ are positive real numbers. It is not hard 
to check that the sequences $\{ a_{n} \}$ and $\{ b_{n} \}$ defined 
recursively by
\ba
a_{n+1} = \frac{a_{n}+b_{n}}{2}, \quad
b_{n+1} = \sqrt{a_{n}b_{n}},
\label{itera1}
\ea
\no
$a_{0} = a,$ and $b_{0}=b$ converge
to a common limit: namely, the {\em arithmetic-geometric mean} of $a$ and 
$b$, denoted by $\text{AGM}(a,b)$. This is a fascinating function; the book
[5] explains its connections with modern algorithms for the
evaluation of $\pi$. 
The reader will find in [1] a survey of 
maps similar to (\ref{itera1}) and an extensive bibliography.

At the turn of the eighteenth century, Gauss [13] was interested in
lemniscates and their lengths. After a numerical calculation, he observed that
\begin{equation}
\frac{1}{\text{AGM}(1, \sqrt{2})} 
\nn
\end{equation}
\no
and 
\begin{equation}
\frac{2}{\pi} \int_{0}^{1} \frac{dx}{\sqrt{1-x^{4}}} 
\nn
\end{equation}
\no
agree to eleven decimal places. (The integral gives the length of a lemniscate.)
With remarkable insight, he discovered that the elliptic integral $G(a,b)$ in
(\ref{elliptic1}) remains invariant if the parameters $(a,b)$ are
replaced with their arithmetic and geometric means; that is,
\ba
G(a,b) & = & G \left( \frac{a+b}{2}, \sqrt{ab} \right).
\label{inv-elliptic}
\ea
\no
Iterating, passing to the limit, and using the invariance of the elliptic 
integral $G$ yields
\ba
G(a,b) & = & \frac{\pi}{2 \, \text{AGM}(a,b)}.
\label{limit-elliptic}
\ea

\no
The convergence of the arithmetic-geometric mean iteration (\ref{itera1}) is
quadratic, meaning that
$|a_{n+1}-\text{AGM}(a,b)| \leq C |a_{n}- \text{AGM}(a,b)|^{2}$ for 
some $C >0$. Thus (\ref{itera1}) leads to a rapid evaluation of the elliptic 
integral $G(a,b)$. This 
iteration has been used for the numerical evaluation of elliptic 
integrals. See [7], [8], [9], [10], or [11] for 
details. The 
algorithm described here could also be used for the numerical
evaluation of rational integrals. 

It was a pleasant surprise when, in the process of analyzing definite integrals
of rational functions, we discovered that
\ba
U_{6} & = & \ift \frac{cx^{4} + dx^{2} + e}{x^{6} + ax^{4} + bx^{2} + 1} \, dx
\nn
\ea
\no
admits a similar invariant transformation. We call this a
{\em rational Landen transformation}. 
In the case of $U_{6}$ the dynamical system (\ref{itera1}) is replaced with
\ba
a_{n+1} & = & \frac{a_{n}b_{n} + 5a_{n} + 5b_{n} + 9}{(a_{n}+b_{n}+2)^{4/3}}, 
\label{itera2} \\
b_{n+1} & = & \frac{a_{n} + b_{n} + 6}{(a_{n}+b_{n}+2)^{2/3}}, \nn
\ea
\no
with similar rules for $c_{n}, \, d_{n},$ and $e_{n}$.  The derivation of 
(\ref{itera2}) appears in [2].

The sequence $(a_{n},b_{n})$ converges to $(3,3)$ 
precisely for those
initial data $(a_{0},b_{0})$ for which the integral 
$U_{6}$ is finite. Moreover, for the 
numerator parameters, we have $(c_{n}, \, d_{n}, \, e_{n}) \to 
(1,2,1)L$, for some real $L$. 
The convergence of this method is discussed 
in [2], [12], and [14]. 
The invariance of $U_{6}$ yields the identity
\ba
U_{6} & = & \frac{\pi}{2 L}
\nn
\ea
\no
exactly as in (\ref{limit-elliptic}). Observe that (\ref{limits2}) is also
of this type: an integral given as the limit of an iterative process. 
Transformations similar to (\ref{itera2}) have been produced in [3]
for any even rational integrand.

Until now all rational Landen transformations were restricted to even 
rational functions. In this paper we 
present the simplest example of a technique
that we expect will extend to the general case (see [17] for details).

The identities (\ref{limits2}) and (\ref{limit-elliptic}) yield 
{\em iterative methods} to evaluate the corresponding integrals. For example,
the first four iterations of the evaluation of
\ba
I & = & \ifft \frac{dx}{4x^{2}+3x +1} \nn
\ea
\no
using (\ref{system1}) are given in Table $1$.

\newpage

\noindent
Table 1. \\

\begin{center}
\begin{tabular}[h]{||c|c|c|c||} 
\hline
$n$ & $a_{n} $ & $b_{n}$ & $c_{n}$ \\ \hline 
0 & 4 & 3 & 1 \\
1 & 1.0731707317 & 0.6585365853 & 1.7317073171 \\
2 & 1.3322738087 & 0.0186646386 & 1.31360991700 \\
3 & 1.3228754233 & 4.644065 $\times 10^{-7} $ & 1.3228758877 \\
4 & 1.3228756555 & 7.154295 $ \times 10^{-21}$  & 1.3228756555 \\
\hline
\end{tabular}
\end{center}

\medskip

The example presented in the table exhibits cubic convergence, faster than the
convergence of the AGM. The exact value of $I$ is $2 \pi/\sqrt{7}$, and  
(\ref{limits2}) yields 
$\lim\limits_{n \to \infty} a_{n}  = \sqrt{7}/2$. The reader can check that 
the value $a_{4}$ gives $\sqrt{7}/2$ 
correct to ten digits of accuracy. 

At the end of the amazing numerical calculation that led him to establish the
invariance for the elliptic integral $G(a,b)$, Gauss commented in his diary
that this {\em will surely open up a whole new field of 
analysis}. This statement is
certainly true. The reader will find in [5] a detailed discussion
of how the arithmetic-geometric mean plays a fundamental role in modern
computations of the digits of $\pi$. This technique 
has also been used in [4]
to create new and efficient methods to evaluate elementary functions. 

Over the years many proofs of (\ref{inv-elliptic}) have been 
discovered. A number
of them can be found in [18]. The authors are particularly fond 
of the succinct proof by D. J. Newman [19]: use
$x = b \, \tan \theta$ and follow with
$x \mapsto x + \sqrt{x^{2} + ab}$. {\em Change of variables is an art}. \\

\medskip

\no
{\bf 3. THE QUADRATIC CASE}. The 
goal of this section is to present the algebraic techniques that 
produce the transformation (\ref{system1}). We scale the 
integrand by multiplying both the numerator and denominator by an 
{\em appropriate} polynomial. This
scaling is one of the main ingredients in the formulation of 
the Landen transformations. The other one will be discussed in the 
next section.

We are motivated by the identities
\ba
U(\tan \theta) = - \frac{\sin (3 \theta)}{\cos^{3}\theta} 
& , & V(\tan \theta)  = - \frac{\cos (3 \theta)}{\cos^{3} 
\theta}, 
\nn
\ea
\no
where 
\ba
U(x) = x^{3} - 3x & ,  & V(x) = 3x^{2}-1.
\nn
\ea

\no
The task is to find coefficients $z_{0}, \, z_{1}, \, z_{2}, \, z_{3}, 
\, z_{4}$ and 
$e_{0}, \, e_{1}, \, e_{2}$ such that
\ba
(ax^{2}+bx+c)(z_{0}x^{4} + z_{1}x^{3} + z_{2}x^{2} + z_{3}x + z_{4})
\label{expansion1}
\ea
\no
can be written as 
\ba
e_{0}U^{2}(x) + e_{1}U(x)V(x) + e_{2}V^{2}(x)
\label{expansion2}
\ea
\no
with unknown coefficients $z_{i}$ and 
$e_{i}$ that are functions of the original parameters 
$a, \, b,$ and $c$. {\em There is so much freedom, it can't be hard}.

Matching (\ref{expansion1}) with (\ref{expansion2}) yields a system of 
seven equations for the eight unknowns. We  use the first five to solve
for the coefficients 
$z_{i}$ in terms of $a, \, b, \, c,$ and the $e_{i}$. To start, comparison
of the constant term in (\ref{expansion1}) and (\ref{expansion2}) gives
\ba
z_{4} & = & c^{-1}e_{2}.
\ea
\no
Using this value, we find that the first-order coefficient  $z_{3}$ 
satisfies $cz_{3}-3e_{1}+bc^{-1}e_{2} = 0$, which yields 
\ba
z_{3} & = & c^{-2}(3ce_{1}-be_{2}).
\ea
\no
The next powers produce
\ba
z_{2} & = & c^{-3} ( 9c^{2}e_{0}-3bce_{1}+b^{2}e_{2} - ace_{2} -6c^{2}e_{2})
\nn
\ea
\no
and 
\ba
z_{1} & = & c^{-4}(-9bc^{2}e_{0}+3b^{2}ce_{1} - 3ac^{2}e_{1}-10c^{3}e_{1} -
b^{3}e_{2} + 2abce_{2} + 6bc^{2}e_{2}), \nn
\ea
\no
respectively. Finally,
\ba
z_{0} & = & c^{-5}( 9b^{2}c^{2}e_{0} 
- 9ac^{3}e_{0}-6c^{4}e_{0} - 3b^{3}ce_{1} + 
6abc^{2}e_{1} + 10bc^{3}e_{1} + b^{4}e_{2} \nn \\
 & & -3ab^{2}ce_{2} + a^{2}c^{2}e_{2} -6b^{2}c^{2}e_{2}+6ac^{3}e_{2} + 
9c^{4}e_{2}). \nn
\ea

This leaves the two equations that arise 
from the two highest powers, which we use to find the parameters 
$e_{i}$. We solve the $x^{5}$ equation for $e_{2}$ 
in terms of the parameters $a, \, b, \, c, \, 
e_{1},$ and $ e_{0}$. Substituting this information into the equation for the
leading term produces
\ba
b( b^{2} - 3(a-c)^{2} )e_{0} & = & a ( 3b^{2} - (a + 3c)^{2}) e_{1}. 
\ea
\no
The system has one degree of freedom, which we exploit to ensure that the 
$z_{i}$ and $e_{i}$ are polynomials in the parameters $a, \, b,$ and $c$.
We initially choose 
$e_{0} = a( (a+3c)^{2} - 3b^{2} )$, from which it follows that
$e_{1}=-b(b^{2}-3(a-c)^{2})$.  This in turn yields
$e_{2}=-c(3b^{2} - (3a+c)^{2})$. 

The expressions for the coefficients $z_{i}$ reduce to the following:
\ba
z_{0} & = & (a+3c)^{2} - 3b^{2} \nn \\
z_{1} & = & 8b(a-3c) \nn \\
z_{2} & = & -6a^{2} + 10b^{2} + 44ac - 6c^{2} \nn \\
z_{3} & = & 8b(c-3a) \nn \\
z_{4} & = & (3a+c)^{2} - 3b^{2} \nn
\ea
\no
and, just to reiterate,
\ba
e_{0} & = & a( (a+3c)^{2} - 3b^{2}) \nn \\
e_{1} & = & b( 3(a-c)^{2} - b^{2}) \nn \\
e_{2} & = & c( (3a+c)^{2} - 3b^{2}). \nn
\ea

\no
In the latter formulas we already see a 
semblance of the iteration (\ref{system1}).

\medskip

\no
{\bf 4. ENTER TRIGONOMETRY}. In 
this section we complete the construction of the Landen transformation
and establish the invariance of the integral (\ref{1.1}) under it. 
We establish the vanishing of a special class of integrals that
appear as intermediate steps in this construction.

We start with (\ref{1.1}) and use the 
change of variables $x = \tan \theta$ to produce 
\ba
I & = & \int_{-\pi/2}^{\pi/2} 
\frac{d \theta}{a \sin^{2} \theta + b \sin \theta \cos \theta + 
c \cos^{2} \theta}. 
\label{trig-integral2}
\ea
\no
The identities
\ba
\tan^{3} \theta - 3 \tan \theta = - \frac{\sin(3 \theta)}{\cos^{3}\theta}, 
\, \, \, 
3 \tan^{2} \theta -1  = - \frac{\cos(3 \theta)}{\cos^{3} 
\theta} \nn
\ea
\no
that were the reason behind the choices for $U$ and $V$ are then used to obtain
\ba
I & = & \sum_{k=0}^{4} z_{4-k} 
\int_{-\pi/2}^{\pi/2} \frac{\sin^{k} \theta \, \cos^{4-k} \theta \, d \theta}
{e_{0} \sin^{2}(3 \theta) + e_{1} \sin( 3 \theta) \cos( 3 \theta) + 
e_{2} \cos^{2}( 3 \theta)}  \label{trig-integral1}
\ea
\no
from the integral (\ref{1.1}) after it has been scaled according to the 
procedure described in section 3.

The elementary identities

\ba
\cos^{4}\theta & = & \tfrac{1}{8} \cos( 4 \theta) + \tfrac{1}{2} 
\cos (2 \theta )
+ \tfrac{3}{8} \label{trig-identities} \\
\cos^{3}\theta \sin \theta & = & \tfrac{1}{8} \sin ( 4 \theta ) + 
\tfrac{1}{4} \sin (2 \theta) \nn \\
\cos^{2} \theta \sin^{2} \theta & = & \tfrac{1}{8} - \tfrac{1}{8} \cos (4 
\theta) \nn \\
\cos \theta \sin^{3}\theta & = & \tfrac{1}{4} \sin (2 \theta) - 
\tfrac{1}{8} \sin (4 \theta) \nn \\
\sin^{4}\theta & = & \tfrac{1}{8} \cos (4 \theta) - \tfrac{1}{2} \cos(2 \theta) 
+ \tfrac{3}{8} \nn
\ea
\no
transform the expression for $I$  to a linear combination of 

\ba
S_{k} & = & \int_{-\pi/2}^{\pi/2} \frac{ \sin( k \theta ) \, d \theta}
{e_{0} \sin^{2} (3 \theta) + e_{1} \sin (3 \theta) \cos (3 \theta) + 
e_{2} \cos^{2} (3 \theta)} \, \, (k = 2, \, 4) \nn
\ea
\no
and 
\ba
C_{k} & = & \int_{-\pi/2}^{\pi/2} \frac{ \cos (k \theta) \, d \theta}
{e_{0} \sin^{2}( 3 \theta) + e_{1} \sin (3 \theta) \cos (3 \theta) + 
e_{2} \cos^{2} (3 \theta)} \, \, (k = 0, \, 2, \, 4) \nn
\ea

\no
The magic of the Landen 
transformations comes from the vanishing of many of these integrals. This 
reduces (\ref{trig-integral1}) to an integral
of the type (\ref{trig-integral2}) with new coefficients, resulting in the
transformation rule (\ref{system1}). Indeed, for even $k$ the  
integrals $S_{k}$ and $C_{k}$ vanish if $k$ is not a multiple of $3$.
To verify this, replace 
$\theta$ with $u = \theta + \pi$ in the definition 
of $S_{k}$. Using 
$\sin(k[u - \pi]) = (-1)^{k} \sin(ku) = \sin (ku)$ and $\cos(k[u - \pi]) = 
(-1)^{k} \cos(ku) = \cos(ku),$ we arrive at 
\ba
S_{k} & = & \int_{\pi/2}^{3\pi/2} \frac{ \sin(ku) \, du}
{e_{0} \sin^{2} (3u)+ e_{1} \sin (3u) \cos (3u) + 
e_{2} \cos^{2} (3u)}. \nn
\ea
\no
Adding this to the original $S_{k}$ and 
taking advantage of the periodicity of the integrand 
we get
\ba
S_{k} & = & \frac{1}{2} \int_{0}^{2 \pi} \frac{ \sin (ku) \, du}
{e_{0} \sin^{2} (3u) + e_{1} \sin (3u) \cos (3u) + 
e_{2} \cos^{2} (3u)}. \nn
\ea
\no
Now, we observe that both $\sin (3u)$ and $\cos (3u)$ are invariant 
under shifts by 
$ 2 \pi/3$ and $4 \pi/3$, so
\ba
6S_{k} & = &  \int_{0}^{2 \pi} \frac{ \sin (ku) + \sin( ku - 2 \pi k/3) +
\sin( ku - 4 \pi k/3)} 
{e_{0} \sin^{2}(3u) + e_{1} \sin (3u) \cos (3u)  + 
e_{2} \cos^{2} (3u)} \, du. \nn
\ea
\no
The numerator in the integrand is the imaginary part of
\ba
e^{iku} + e^{i(ku- 2 \pi k/3)} + e^{i(ku - 4 \pi k/3)} & = & 
e^{i k u} \left( 1  + e^{- 2 \pi k i /3} + e^{- 4 \pi k i /3} \right), \nn
\ea
\no
and the last sum is $3$ or $0$ depending on whether $3$ divides $k$ or not.

\medskip

We conclude that the only terms that contribute to (\ref{trig-integral1}) are
the constants in (\ref{trig-identities}). Therefore
\ba
I & = & \frac{1}{16} \int_{0}^{2 \pi} \frac{3z_{4} + z_{2} + 3z_{0}}
{e_{0} \sin^{2} (3u) + e_{1} \sin (3u) \cos (3u) + 
e_{2} \cos^{2} (3u)} \, du, \nn
\ea
\no
where we have again appealed to periodicity to extend the 
integral to $[0, 2 \pi]$. The 
change of variables $\theta = 3u$ leads to 
\ba
I & = & \frac{1}{8} \int_{-\pi/2}^{\pi/2} \frac{3z_{4} + z_{2} + 3z_{0}}
{e_{0} \sin^{2} \theta  + e_{1} \sin \theta  \cos \theta  + 
e_{2} \cos^{2} \theta}  \, d \theta, \nn
\ea
\no
so we have returned to the original form
(\ref{trig-integral2}) but with different coefficients.  The result in 
(\ref{invariance1}) is obtained by using 
$x = \tan \theta$ and the following identities:
\ba
\frac{8e_{0}}{3z_{4}+z_{2}+3z_{0}} &  = & 
a \left( \frac{(3a+c)^{2}-3b^{2}}{(3a+c)(a+3c)-b^{2}} \right) \label{iden0} \\
\frac{8e_{1}}{3z_{4}+z_{2}+3z_{0}} &  = & 
b \left( \frac{3(a-c)^{2} - b^{2}}{(3a+c)(a+3c)-b^{2}} \right) \nn \\
\frac{8e_{2}}{3z_{4}+z_{2}+3z_{0}} &  = & 
c \left( \frac{(a+3c)^{2}-3b^{2}}{(3a+c)(a+3c)-b^{2}} \right). \nn 
\ea

\medskip

\no
{\bf 5. THE ANALYSIS OF CONVERGENCE}. In 
the last two sections we have shown the invariance of (\ref{1.1}) 
under the Landen transformation (\ref{system1}). We now conclude by 
establishing the convergence of its iterates as in (\ref{limits1}). In 
particular, we show that the error
\ba
e_{n} & := & ( a_{n} - \tfrac{1}{2}\sqrt{4ac - b^{2}}, \, b_{n}, \, 
c_{n} - \tfrac{1}{2}\sqrt{4ac - b^{2}} ) 
\nn
\ea
\no
satisfies $e_{n} \to 0$ as $n \to \infty$. Moreover, we demonstrate 
cubic convergence:
\ba
\| e_{n+1} \| & \leq & C \| e_{n} \|^{3}
\label{cubic-error}
\ea
\no
for some positive constant  $C$. 

The analysis of convergence is simpler in the variables 
$x = a+c, \, y=b$, and $z=a-c$. The dynamical system (\ref{system1})
translates to 
\ba
x_{n+1} & = & x_{n} \left[ \frac{4x_{n}^{2} -3z_{n}^{2} - 3y_{n}^{2}}
{4x_{n}^{2} -y_{n}^{2} - z_{n}^{2}} \right], \label{system2} \\
z_{n+1} & = & z_{n} \left[ \frac{z_{n}^{2} -3y_{n}^{2}}
{4x_{n}^{2} -y_{n}^{2} - z_{n}^{2}} \right], \nn  \\
y_{n+1} & = & y_{n} \left[ \frac{3z_{n}^{2} -y_{n}^{2}}
{4x_{n}^{2} -y_{n}^{2} - z_{n}^{2}} \right], \nn  
\ea
\no
with initial conditions $x_{0} =x, \, y_{0}=y$, and $z_{0}=z$. 

We now prove that
\ba
\lim\limits_{n \to \infty} x_{n} = \sqrt{x^{2}-y^{2}-z^{2}}, \quad
\lim\limits_{n \to \infty} y_{n} = 
\lim\limits_{n \to \infty} z_{n} = 0 
\nn
\ea
\no
or, equivalently, that 
\ba
\lim\limits_{n \to \infty} \left( x_{n} - \sqrt{x^{2}-y^{2}-z^{2}} \right)^{2} 
+ y_{n}^{2} + z_{n}^{2} & = & 0.
\label{lastlimit}
\ea
\no
This is equivalent to (\ref{limits1}), so it will finish the proof of 
convergence.

To complete the change of variables we use the 
invariance of the discriminant (\ref{discriminant1}) to obtain
\ba
x_{n}^{2} - y_{n}^{2} - z_{n}^{2} = x^{2}-y^{2}-z^{2} = 4ac - b^{2},
\nn
\ea
\no
and we write $w = \sqrt{4ac-b^{2}}$. 
The first equation of iteration (\ref{system2}) becomes
\ba
x_{n+1} & = &  x_{n} 
\left[ \frac{x_{n}^{2} + 3w^{2}}{3x_{n}^{2} + w^{2}} \right], \label{system3}
\ea
\no
with initial conditions $x_{0} = a+c \, (> 0)$. The required limit in 
(\ref{lastlimit}) is now
\ba
\lim\limits_{n \to \infty} x_{n}(x_{n}- w) & = & 0,
\label{last-limit}
\ea
\no
with $y_{n}$ and $z_{n}$ absorbed into the constant $w$. 
{\em The number of variables has been reduced from three to one}. 

Using one last change of variables, 
$q_{n} = -ix_{n}/w$, we reduce (\ref{system3}) to

\ba
q_{n+1} & = & \frac{q_{n}^{3}-3q_{n}}{3q_{n}^{2}-1} = 
\frac{U(q_{n})}{V(q_{n})}. \label{system4}
\ea
\no
What we need to prove in order to establish (\ref{last-limit}) is that 
$q_{n} \to -i$. (The 
polynomials $U$ and $V$ introduced in section 3 have 
miraculously reappeared!) The trigonometric identity 
\ba
\frac{U( \cot \theta)}{V( \cot \theta )} & = & \cot( 3 \theta), 
\nn
\ea
\no
coupled with a representation of the initial condition as
\ba
q_{0} & = & \cot( i t )
\label{qzero}
\ea
\no
for some $t \, (0 < t < \infty)$, shows that (\ref{system4}) simplifies to 
\ba
q_{1}  & = & \frac{U( \cot \, it )}{V( \cot \, it )} = \cot( 3 i t)
\nn
\ea
\no
and, in general,
\ba
q_{n} = \text{cot}(3^{n}it) = -i \frac{e^{2t \, 3^{n}} +1}{e^{2t \, 3^{n}} -1}.
\nn
\ea
\no
We conclude that $q_{n} \to -i$, whence $x_{n} \to w$, as desired. 

To verify (\ref{qzero}), we write $q_{0} = -id$ with
$d = (a+c)/\sqrt{4ac-b^{2}}$. Now recall that $4ac-b^{2} > 0.$ An elementary
argument shows that 
$d \geq 1$, so we can take 
\ba
t & = & \text{coth}^{-1}(d) = \frac{1}{2} \ln \frac{d+1}{d-1}. 
\nn
\ea

The fact that the convergence is cubic follows directly from 
\ba
| q_{n} + i | & = & \frac{2}{e^{2 t 3^{n}}-1}, 
\nn
\ea
\no 
which decreases to $0$ like $e^{-2t3^{n}}$. This implies (\ref{cubic-error}) 
and completes the proof of convergence.

\bigskip

\no
{\bf ACKNOWLEDGMENTS}. The authors wish to thank the referees for a careful 
reading of the original manuscript. The
second author acknowledges the partial support of
NSF award DMS-0409968. The first author was partially 
supported as a 
graduate student by the same grant.

\bigskip

\no
REFERENCES \\

\no
1. J. Arazy, T. Claesson, S. Janson, and J. Peetre, Means and their iterations, in {\em Proceedings of the Nineteenth Nordic Congress of Mathematicians, 
Reykjavik}, Icelandic Mathematical Society, 1984, pp. 191-212.  \\

\no
2. G. Boros and V. Moll, A rational Landen transformation. The case of degree
$6$, in {\em Analysis, Geomtery, Number Theory: The Mathematics of Leon 
Ehrenpreis}, Contemporary Mathematics, vol. 251, M. Knopp, G. Mendoza, 
E. T. Quinto, E. L. Grinberg, and S. Berhanu, eds., American Mathematical 
Society, Providence, 2000, pp. 83-89. \\

\no
3. $\quad \quad $, Landen transformations and the integration of 
rational functions, {\em Math. Comp.} {\bf 71} (2001)  649-668. \\

\no
4. J. M. Borwein and P. B. Borwein, The arithmetic-geometric mean and fast 
computation of elementary functions, {\em SIAM Review} {\bf 26} (1984) 
351-366. \\

\no
5. $\quad \quad$, {\em Pi and the AGM - A study in analytic 
number theory and computational complexity}, Wiley, New York, 1987. \\

\no
6. R. Burlisch, Numerical calculation of elliptic integrals and functions, 
{\em Numer. Math.} {\bf 7} (1965) 78-90. \\

\no
7. $\quad \quad$, Numerical calculation of elliptic 
integrals and functions, II, 
{\em Numer. Math.} {\bf 7} (1965) 353-354. \\

\no
8. $\quad \quad$, Numerical 
calculation of elliptic integrals and functions, III, 
{\em Numer. Math.} {\bf 13} (1969) 305-315. \\

\no
9. B. C. Carlson, Algorithms involving arithmetic and geometric means, {\em Amer. Math. Monthly} 
{\bf 78} (1971) 496-505. \\

\no
10. $\quad \quad$, Computing elliptic integrals by duplication, {\em Numer. 
Math.} {\bf 33} (1979) 1-16. \\

\no
11. $\quad \quad$, Numerical computation of real or complex elliptic integrals,
{\em Numer. Algorithms} {\bf 10} (1995) 13-26. \\

\no
12. M. Chamberland and V. Moll, Dynamics of the degree six Landen 
transformation, {\em Discrete and Continuous Dynamical Systems} {\bf 15} (2006) 21-37. \\

\no
13. K. F. Gauss, Arithmetische Geometrische Mittel in {\em Werke}, 
vol. 3 (1799) 361-432; reprinted by Olms, Hildescheim, 1981. \\

\no
14. J. Hubbard and V. Moll, A geometric view of rational Landen 
transformations, {\em Bull. London Math. Soc.} {\bf 35} (2003) 293-301. \\

\no
15. C. G. J. Jacobi, Fundamenta nova theoriae funcionum ellipticarum in 
{\em Gesammelte Werke}, vol. I (1829) 49-239; reprinted by 
Chelsea Publishing Company, New York, 1969. \\

\no
16. F. Klein, {\em Developments of Mathematics in the $19^{th}$ Century}. 
Reprinted by Trans. Math. Sci. Press, R. Hermann ed., Brookline, MA, 1979. \\

\no
17. D. Manna, {\em Landen transformations}. Ph.D. Thesis, Tulane University, 
2006. \\

\no
18. H. McKean and V. Moll, {\em Elliptic Curves: Function Theory, Geometry, 
Arithmetic}. Cambridge University Press, New York, 1997. \\

\no
19. D. J. Newman, A simplified version of the fast algorithm of Brent and
Salamin, {\em Math. Comp.} {\bf 44} (1985) 207-210.

\bigskip
\no
{\em Department of Mathematics, Tulane University, New Orleans, LA 70118} \\

\no
{\em dmanna@math.tulane.edu} and {\em vhm@math.tulane.edu}

\end{document}